\documentclass{svproc}
\usepackage{url}
\usepackage{amsmath}
\usepackage{amsfonts}
\usepackage{graphicx}
\usepackage{xcolor, ulem}
\usepackage{caption}
\usepackage{subcaption}

\newcommand*{\OO}[1]{\mathbf{O}_{#1}}                   
\newcommand*{\Gr}[1]{\mathbf{Gr}_{#1}}        
\newcommand*{\St}[1]{\mathbf{St}_{#1}}        

\newcommand{\R}{\mathbb{R}}

\def\tr{{\rm{tr }}}

\newcommand*{\so}[1]{\mathfrak{so}(#1)}					
\newcommand*{\soP}[1]{\mathfrak{so}_P(#1)}					

\begin{document}
\mainmatter              
\title{Comparison of two numerical methods for Riemannian cubic polynomials on Stiefel manifolds}
\titlerunning{Numerical methods for Riemannian cubic polynomials}  
%
\author{Alexandre Anahory Simoes\inst{1} \and Leonardo Colombo\inst{2} \and F\'atima Silva Leite\inst{3}}
\authorrunning{Anahory et al.} 
%
\tocauthor{Leonardo Colombo, F\'atima Silva Leite, Alexandre Anahory Simoes}
\institute{IE University, Madrid, 28027, Spain,\\
\email{alexandre.anahory@ie.edu}
\and
Center for Automation and Robotics (UPM-CSIC), Spain\\ \email{leonardo.colombo@csic.es}
\and
Institute of Systems and Robotics, and 
Department of Mathematics, University of Coimbra, Portugal\\ \email{fleite@mat.uc.pt}
}

\maketitle              

\begin{abstract}
In this paper we compare two numerical methods to integrate Riemannian cubic polynomials on the Stiefel manifold $\St{n,k}$.  The first one is the adjusted de Casteljau algorithm, and the second one is a symplectic integrator constructed through discretization maps. In particular, we choose the cases of  $n=3$ together with $k=1$ and $k=2$. The first case is diffeomorphic to the sphere and the quasi-geodesics appearing in the adjusted de Casteljau algorithm are actually geodesics. The second case is an example where we have a pure quasi-geodesic different from a geodesic. We provide a numerical comparison of both methods and discuss the obtained results to highlight the benefits of each method.
\end{abstract}
\section{Introduction}

Generating polynomial curves and polynomial splines on manifolds was motivated by problems related to path     planning of certain mechanical systems, such as spacecraft and underwater vehicles, whose configuration spaces are non-Euclidean manifolds {\cite{Park:1995aa}}, \cite{camarinha2001geometry}, \cite{bloch2021dynamic}, {but the impact of the results quickly expanded to other areas of technology}. A Riemannian cubic polynomial is a smooth curve on a Riemannian manifold $Q$ that minimizes the cost functional
$$\mathcal{J}= \int_{0}^{1}  \left\langle \frac{D \dot{c}}{dt},\frac{D \dot{c}}{dt} \right\rangle  \ dt,$$
where $\langle \cdot, \cdot \rangle$ denotes the Riemannian metric and $\frac{D \dot{c}}{dt}$ denotes the covariant acceleration of the curve $c$ with respect to the Levi-Civita connection.  In addition, it satisfies a two-point
boundary value problem
(initial and final points and velocities are prescribed). 
Without loss of generality, we parameterize the curves over the interval $[0,1]$. {It turns out that the corresponding Euler-Lagrange equations are hard to solve and to overcome such difficulty other methods to find cubic polynomials on manifolds have been proposed, in particular geometric algorithms that generate such curves  from four distinct points
 $x_0,x_1, x_2, x_3$ on the configuration manifold, the first and last being respectively
 the initial and final prescribed points and the other two are only auxiliary points for the  geometric construction, but  are also related to the specified velocities}.

The Stiefel manifold, {consisting of all $n\times k$ ($1\leq k\leq n$) matrices with orthonormal columns}, is important in various applications, including computer vision \cite{turaga2008statistical}, neural networks
\cite{nishimori2005learning}, and statistics \cite{chakraborty2019statistics}. {For instance, its relevance in pattern
recognition is due to the fact that features and patterns that describe visual objects may be represented as elements in those non-Euclidean spaces}. These geometric representations facilitate the analysis of the underlying geometry of the data. Stiefel manifolds are related to Grassmannian manifolds, {consisting of all k-dimensional subspaces in $\mathbb{R}^{n}$, since a point in the Stiefel manifold identifies exactly what frame (basis of orthonormal vectors) is used to specify a particular subspace}.

In this paper, we will provide a comparison between a recently proposed {geometric} method to approximate Riemannian cubic polynomials on a Stiefel manifold with a symplectic integrator for the same {variational problem}. The {approach} we propose consists of a symplectic integrator generated from a choice of a retraction map \cite{Barbero-Linan:2023aa}, \cite{simoes2023higher}. This formulation seems to be {well} suited for a wide range of problems in geometric control theory, namely optimal control problems involving higher-order derivatives whose state space is a Riemannian space with additional structure such as that of a Lie group or that of a homogeneous space.

 In particular, we choose the cases of  $n=3$ together with $k=1$ and $k=2$. The first case is diffeomorphic to the sphere and the quasi-geodesics {defined} in \cite{Krakowski:2017aa} are actually geodesics. The second case is an example where we have a pure quasi-geodesic different from a geodesic. We provide a numerical comparison of both methods and discuss the obtained results to highlight the benefits of each method.

\section{The Stiefel manifold}\label{sec-stiefel}

\subsection{Background \& notations}

The Stiefel manifold of orthonormal $k$-frames in $\R^n$ has the following matrix representation:
\begin{equation}
	\St{n,k}=\{S\in\R^{n\times k} \mid S^\top S=I_k\} .
\end{equation}
\

We note that, if $S\in \St{n,k}$, then $P=SS^\top$ is a point in the Grassmann manifold $\Gr{n,k}$ consisting of all $k$-dimensional subspaces of $\R^n$, when its matrix representation   $\Gr{n,k}=\{P\in\R^{n\times n} \mid P=P^\top, \, P^2=P, \, rank(P)=k\} $ is considered.
 
In what follows, $\so{n}$ denotes the set of $n\times n$ skew-symmetric matrices, and for $P\in \Gr{n,k}$, the notation $\soP{n}$ is used for the vector subspace of $\so{n}$ defined by
$$\soP{n} =\{ X\in  \so{n}\mid XP + PX =X \}.$$ 
We also note that 
$$
	T_P\Gr{n,k} = \{ [X,P]\mid X\in \soP{n}\} ,
$$
where $[\cdot,\cdot]$ denotes the commutator of matrices.

The tangent space to the Stiefel manifold at a point  $S\in \St{n,k}$ is given by
\begin{equation}\label{eq-stiefel-tangent}
	T_S\St{n,k} = \{V\in\R^{n\times k}\mid V^\top S+S^\top V=0\} ,
\end{equation}
but another useful representation of the tangent space appeared in ~\cite[Proposition 5]{Krakowski:2017aa} and is recalled here. Suppose  $S\in \St{n,k}$ and  $P:= SS^\top \in \Gr{n,k}$. Then,
\begin{equation}\label{new-param}
T_S\St{n,k}=\{ XS+S\Omega \mid X\in \soP{n}, \,\,\,   \Omega\in \so{k}\}.
\end{equation}
Moreover, if $V=XS+S\Omega \in T_S\St{n,k}$, then 
\begin{equation}\label{new-param1}
X= VS^\top -SV^\top +2SV^\top SS^\top , \quad   \Omega= S^\top V. 
 \end{equation}

We will consider the Stiefel manifold equipped with the metric defined by
\begin{equation}\label{innerproduct}
	\langle V_1, V_2\rangle = \tr \big( V_1^\top (I_n-\tfrac{1}{2}SS^\top ) V_2\big),\quad V_1,V_2\in T_S\St{n,k}.
\end{equation}
In some literature, for instance in~\cite{Edelman:98}, this is called the \emph{canonical metric}. Also, among the family of $\alpha$-metrics studied in~\cite{Hueper:21}, this  corresponds to the value $\alpha =0$. 


The geometric de Casteljau algorithm to generate polynomials on Riemannian manifolds is based on successive geodesic interpolation. {To deal with situations for which an explicit expression for a geodesic  that joins two points} is not available, the authors of \cite{Krakowski:2017aa} used quasi-geodesics, instead of geodesics,  to modify the de Casteljau algorithm to generate quadratic polynomials and splines on Stiefel manifolds. In \cite{preprintM}, the adjusted de Casteljau algorithm is used to approximate Riemannian cubic polynomials and solve {related} interpolation problems numerically.

\subsection{Retractions and quasi-geodesics on Stiefel manifolds}\label{sec-retractions}

{We first recall from~\cite{AbMaSeBookRetraction}} the notion of retraction map on a Riemannian manifold.
\begin{definition}\label{retraction}
A retraction $R$ on the Stiefel manifold $\St{n,k}$ is a smooth mapping  from the tangent bundle $T\St{n,k}$ to $\St{n,k}$ that, when restricted to each tangent space at a point $S\in \St{n,k}$ (restriction denoted  by $R_S$),  satisfies the following properties:
\begin{enumerate}
\item[(i)] $R_S(0)=S$;
\item[(ii)] $dR_S(0)= I$, where $dR_S(0)$ stands for the tangent map of $R_S$ at $0\in T_S\St{n,k}$.
\end{enumerate}
\end{definition}
If $V \in T_S\St{n,k}$, one can define a smooth curve $\beta _V\colon t\mapsto R_S(tV)$ associated to the retraction $R$. The curve $\beta_V$ which satisfies $\beta_V(0)=S$ and $\dot \beta_V(0)=V$ is called a \emph{quasi-geodesic}. 
Next, we  define a particular retraction and corresponding  quasi-geodesic on the Stiefel manifold, and  list some of their interesting properties. 
Proofs and more details can be found in~\cite{Krakowski:2017aa}. We may use indistinctly $\exp{(A)}$ or $e^A$ to denote the  matrix exponential of a matrix $A$ and $log$ will denote the principal matrix logarithm.

If $S$, $X$ and $\Omega$ are as in equation \eqref{new-param1}, then the mapping $R\colon T\St{n,k}\to \St{n,k}$ whose restriction to $T_S\St{n,k}$ is defined by 
$R_S(V)=e^XSe^\Omega$ is a retraction on the Stiefel manifold, {with associated quasi-geodesic  $\beta: [0,1] \rightarrow \St{n,k}$, $t\mapsto e^{tX} S e^{t\Omega}$, which  satisfies}
$$\beta(0) = S;\quad \dot\beta(t) = e^{tX} (XS + S\Omega) e^{t\Omega}; \quad \ddot\beta(t) = e^{tX} (X^2S + 2XS\Omega + S\Omega^2) e^{t\Omega}.
 $$


In the next proposition, adapted from \cite{preprintM}, the initial velocity of a quasi-geodesic that joins two points {$S_0$ and $S_1$ is explicitly written in terms of these endpoints}.
\begin{proposition} \label{theo-quasi}
Let $S_0$ and $S_1$ be two distinct points in $\St{n,k}$ so that, for $i=0,1$, $P_i=S_iS_i^\top \in \Gr{n,k}$. Then, if
\begin{equation}\label{eq-q-geodesic} X = \frac{1}{2}\log \big((I-2S_1S_1^\top )(I-2S_0S_0^\top)\big) \,\,\mbox{ and }\,\, 
		\Omega = \log \big(S_0^\top e^{-X} S_1\big), \end{equation}
the quasi-geodesic $\beta: [0,1]\mapsto \St{n,k}$ defined by
\begin{equation}\label{eq-q-geodesic1}
	\beta(t) = e^{tX} S_0 e^{t\Omega}
\end{equation}
is a smooth curve connecting $S_0$ to $S_1$.
\end{proposition}
\begin{remark} 
Note that the matrices $X$ and $\Omega$ in \eqref{eq-q-geodesic} are only well defined if the logarithm exists. This is always  guaranteed if the points $S_0$ and $S_1$ are sufficiently close. In addition, the quasi-geodesic  defined above is a true geodesic w.r.t. the metric \eqref{innerproduct} only if $X=0$ or $\Omega=0$. In particular, these two situations occur when $k=1$ or $k=n$, {in which cases $\St{n,1}$ is the unit sphere $S^n$  or the orthogonal group $\OO{n}$.}
\end{remark}

  \section{The adjusted de Casteljau Algorithm}
\label{sec:de-cast-algor-1}
%

We briefly describe the adjusted de Casteljau Algorithm to generate
cubic polynomials on Riemannian manifolds, assuming that they are geodesically complete \cite{Crouch:1999ww}. The next algorithm results from replacing geodesics by quasi-geodesics in the generalized de Casteljau algorithm {concisely} described in \cite{H-SL-23} and adapting it to the construction of cubic polynomials as in \cite{preprintM}.

\begin{problem} \label{problem} Find a smooth curve $\gamma: [0,1]\rightarrow \St{n,k}$ satisfying the following boundary conditions:
\begin{equation}
\gamma (0)=S_0, \quad \gamma (1)=S_3, \quad \dot\gamma (0)=V_0, \quad \dot\gamma (1)=V_3,
\end{equation}
where $S_0,S_3$ are given points in $\St{n,k}$, and $V_0\in T_{S_0}\St{n,k}$ and  $V_3\in T_{S_3}\St{n,k}$ are given tangent vectors.
\end{problem}

The adjusted de Casteljau algorithm is used to generate a curve that solves {the  previous problem. Following \cite{preprintM},  we first define matrices $X_i$ and $\Omega_i$, with  $i$ an integer running from $0$ to $5$, not necessarily in this particular order. In the way,  we also find the control points $S_1, S_2$ from the given data.} 
\begin{equation}\label{X0-Omega0}
X_0= \frac{1}{3}\big( V_0S_0^\top -S_0V_0^\top +2S_0V_0^\top S_0S_0^\top \big), \qquad \Omega_0= \frac{1}{3}S_0^\top V_0 .
\end{equation}
So, $S_1=e^{X_0}S_0e^{\Omega_0}$ defines the \emph{first control point}. Using

 
\begin{equation}\label{X2-Omega2}
X_2= -\frac{1}{3}\big( V_3S_3^\top -S_3V_3^\top +2S_3V_3^\top S_3S_3^\top \big), \qquad \Omega_2= -\frac{1}{3}S_3^\top V_3,
\end{equation}
the point $S_2=e^{X_2}S_3e^{\Omega_2}$ is the \emph{second control point} which enables to define 
 

\begin{equation*}\label{X1-Omega1}
X_1= \frac{1}{2}\log\big( (I-2S_2S_2^\top )(I-2S_1S_1^\top )  \big), \qquad \Omega_1= \log\big( S_1^\top e^{-X_1}S_2 \big) .
\end{equation*}
Finally, define
\begin{equation*}
\begin{array}{l}
X_3(t)= \frac{1}{2}\log\big( (I-2e^{tX_1}S_1 S_1^\top e^{-tX_1}) (I-2e^{tX_0}S_0 S_0^\top e^{-tX_0}) \big) ,\\
\\
\Omega_3(t)= \log\big( e^{-t\Omega_0}S_0^\top e^{-tX_0}e^{-X_3(t)} e^{tX_1}S_1  e^{t\Omega_1} \big) ,\\
\\
X_4(t)= \frac{1}{2}\log\big( (I-2e^{-tX_2}S_2 S_2^\top e^{tX_2}) (I-2e^{tX_1}S_1 S_1^\top e^{-tX_1}) \big) ,\\
\\
\Omega_4(t)= \log\big( e^{-t\Omega_1}S_1^\top e^{-tX_1}e^{-X_4(t)} e^{-tX_2}S_2  e^{-t\Omega_2} \big) ,\\
\\
X_5(t)= \frac{1}{2}\log\big( (I-2e^{tX_4(t)}e^{tX_1}S_1 S_1^\top e^{-tX_1}e^{-tX_4(t)})\\
\hspace{6 cm} (I-2e^{tX_3(t)}e^{tX_0}S_0 S_0^\top e^{-tX_0}e^{-tX_3(t)}) \big) ,\\
\\
\Omega_5(t)= \log\big( e^{-t\Omega_3(t)}e^{-t\Omega_0}S_0^\top e^{-tX_0}e^{-tX_3(t)}e^{-X_5(t)}e^{tX_4(t)}e^{tX_1}S_1  e^{t\Omega_1} e^{t\Omega_4(t)}\big).
\end{array}
\end{equation*}    

\begin{theorem} \label{theo}
Let $X_i$ and $\Omega_i$ be given by the last formulas. Then, the geometric cubic polynomial  $\gamma: [0,1]\rightarrow \St{n,k}$  defined by
\begin{equation}\label{cubic}
\begin{array}{lcl}
 \gamma (t) &= &e^{tX_5(t)} e^{tX_3(t)}e^{tX_0}S_0e^{t\Omega_0}e^{t\Omega_3(t)} e^{t\Omega_5(t)},
 \end{array}
 \end{equation}
 solves \bf{Problem \ref{problem}}.
\end{theorem}
{The  proof of this theorem can be found in preprint \cite{preprintM}.}

%
\section{Symplectic integrators arising from Discretization maps}
%


Next, we define a generalization of the retraction map that allows a discretization of the tangent bundle of the configuration manifold leading to the construction of numerical integrators as described in~\cite{Barbero-Linan:2023aa}. Given a point and a velocity, we obtain two nearby points that are not necessarily equal to the initial base point. 

\begin{definition} \label{def:DiscreteMap2} Let $Q$ be a smooth manifold. A map 
	$R_d\colon U\subset TQ\rightarrow Q\times Q$ given by \begin{equation*}
	R_d(q,v)=(R^1(q,v),R^2(q,v)),
	\end{equation*} 
	where  $U$ is an open neighborhood of the zero section $0_q$ of $TQ$, 
	defines a {\it discretization map on $Q$} if it satisfies 
	\begin{enumerate}
		\item $R_d(q,0)=(q,q)$,
		\item $T_{0_q}R^2_q-T_{0_q}R^1_q\colon T_{0_q}T_qQ\simeq T_qQ\rightarrow T_qQ$ is equal to the identity map on $T_qQ$ for any $q$ in $Q$, where $R^a_q$ denotes the restrictions of $R^a$, $a=1,2$, to $T_qQ$.
	\end{enumerate}
\end{definition}
Thus, the discretization map $R_d$ is a local diffeomorphism from some neighborhood of the zero section of $TQ$. If $R^1(q,v)=q$, the two properties in Definition~\ref{def:DiscreteMap2} guarantee that \sout{the} both properties in Definition~\ref{retraction} are satisfied by $R^2$. 

The two following examples will be used later during the numerical integration of Riemannian cubic polynomials.

\begin{example}\label{intial:point:discretization}
    The initial point rule on an Euclidean vector space is given by the discretization map:
    $R_d(q,v)=\left( q, q+v\right).$
\end{example}

\begin{example}\label{mid:point:discretization} The mid-point rule on an Euclidean vector space can be recovered from the following discretization map:
$R_d(q,v)=\left( q-\dfrac{v}{2}, q+\dfrac{v}{2}\right).$
\end{example}

Different discretizations might be used to produce different symplectic discretizations of Hamiltonian flows \cite{Barbero-Linan:2023aa}. We briefly outline the construction of these numerical methods. As the Hamiltonian vector field takes values on $TT^*Q$, the discretization map must be on $T^*Q$. If one chooses the cotangent lift of the discretization map $R_d\colon TQ\rightarrow Q\times Q$ composed with additional canonical isomorphisms to place it in the correct spaces just as defined in \cite{Barbero-Linan:2023aa}, we obtain a map denoted by
$R^{T^*}_d: TT^*Q \rightarrow T^*Q\times T^*Q$,  which is not only a discretization but also a symplectomorphism between $(T(T^*Q), {\rm d}_T \omega_Q)$ and $(T^*Q\times T^*Q, \Omega_{12})$, where $\omega_{Q}$ is the canonical symplectic {form} on $T^{*}Q$, ${\rm d}_T \omega_Q$ is the symplectic form on $TT^{*}Q$ obtained from lifting $\omega_{Q}$, and $\Omega_{12}$ is the canonical symplectic structure on the product of two symplectic manifolds. These additional properties allow to prove the following result defining a symplectic numerical scheme:
\begin{proposition}\label{Prop:symplectic:integrator}
		If $R_{d}$ is a discretization map on $Q$ and $H:T^{*}Q \rightarrow \R$ is a Hamiltonian function, then the equation
        \begin{align*}
            &(R_{d}^{T^*})^{-1}(q_{0},p_{0},q_{1},p_{1}) = h X_{H} \left[ \tau_{T^{*}Q} \circ (R_{d}^{T^*})^{-1}(q_{0},p_{0},q_{1},p_{1}) \right]
        \end{align*}
        written for the cotangent lift of $R_{d}$ is a symplectic integrator, where $X_{H}$ is the associated Hamiltonian vector  field on $T^{*}Q$.
	\end{proposition}

\begin{example}\label{example3} On $Q={\mathbb R}^n$ the discretization map 
	$R_d(q,v)=\left(q-\frac{1}{2}v, q+\frac{1}{2}v\right)$ is cotangently lifted to
		$$R_d^{T^*}(q,p,\dot{q},\dot{p})=\left( q-\dfrac{1}{2}\,\dot{q}, p-\dfrac{\dot{p}}{2}; \; q+\dfrac{1}{2}\, \dot{q}, p+\dfrac{\dot{p}}{2}\right)\, .$$
  The initial point discretization map $R_d(q,v)=\left(q, q+v\right)$ is cotangently lifted to
		$$R_d^{T^*}(q,p,\dot{q},\dot{p})=\left( q, p-\dot{p}, \; q+\dot{q}, p\right)\, .$$
\end{example}

The previous proposition adapts perfectly to our case since {Riemannian} cubic polynomials can be seen as the projection of an Hamiltonian flow on $T^{*}(T\St{n,k})$ (see, e.g. \cite{CROUCH200013}). As a result, we may use the previous {proposition} to construct a symplectic integrator for the Hamiltonian version of Riemannian cubic polynomials.

 \subsection{Hamiltonian in the Stiefel manifold}

 {Getting} a general expression for the Hamiltonian function associated with Riemannian cubic polynomials might be a hard task. To reduce the difficulty of the problem and give a comparison between different methods, we will focus on just two cases: $n=3$ together with $k=1$ and $k=2$. The first case is diffeomorphic to the sphere and the quasi-geodesics appearing in the adjusted de Casteljau algorithm are actually geodesics. The second case is an example where we have a pure quasi-geodesic different from a geodesic. Next, we will write the Hamiltonian appearing in Proposition \ref{Prop:symplectic:integrator} in a coordinate chart.

The Stiefel manifold $\St{3,1}$ is diffeomorhpic to the sphere. Thus, we may choose as a chart the parametrization of the sphere using spherical coordinates
\begin{equation*}
    \begin{split}
        \psi: ]0,\pi[ \times ]0, 2\pi[ & \rightarrow S^{2} \\
        (\theta, \phi) & \mapsto (\cos \phi \sin \theta, \sin \phi \sin \theta, \cos \theta)
    \end{split}
\end{equation*}

Riemannian cubic polynomials on the sphere are the projection to $S^{2}$ of the Hamiltonian flow with respect to the Hamiltonian function $H:T^{*}(TS^{2}) \rightarrow \R$ given in coordinates by:
$$H(\theta, \phi, \dot{\theta}, \dot{\phi}, p_{\theta}, p_{\phi}, p_{\dot{\theta}}, p_{\dot{\phi}}) = \frac{1}{2}\dot{\phi}^2 p_{\dot{\theta}} \sin(2 \theta) + \dot{\phi} p_{\phi} + \dot{\theta} p_{\theta} + \frac{1}{2} p_{\dot{\theta}}^2 + \frac{-\dot{\phi} \dot{\theta} p_{\dot{\phi}} \sin(2 \theta) + \frac{1}{2} p_{\dot{\phi}}^2}{\sin^{2}(\theta)}.$$

In the second case, to choose a chart on $\St{3,2}$ we should observe that this manifold is diffeomorphic to the unitary bundle of the sphere, i.e., the set of points in $TS^{2}$ satisfying $\|v_{q}\|=1$ for each $v_{q}\in T_{q}S^{2}$ and the norm is the euclidean norm obtained from the canonical inclusion of $S^{2}$ on $\R^{3}$.

Under this identification, we proceed with the following choice of chart on $\St{3,2}$
\begin{equation*}
    \begin{split}
        \Psi: ]0,\pi[ \times ]0, 2\pi[ \times ]0, 2\pi[ & \rightarrow \St{3,2} \\
        (\theta, \phi, \psi) & \mapsto A(\theta, \phi, \psi)
    \end{split}
\end{equation*}
with
$$A(\theta, \phi, \psi) = \begin{bmatrix}
    \cos \phi \sin \theta &  -\sin(\phi) \sin(\psi) + \cos(\phi) \cos(\psi) \cos(\theta) \\
    \sin \phi \sin \theta &
    \sin(\phi) \cos(\psi) \cos(\theta) + \sin(\psi) \cos(\phi) \\
    \cos \theta & -\sin(\theta) \cos(\psi)
\end{bmatrix}$$

Riemannian cubic polynomials are approximated as the projection to $\St{3,2}$ of the Hamiltonian flow on $T^{*}(T\St{3,2})$ with respect to the Hamiltonian function
\begin{tiny}
\begin{equation*}
    \begin{split}
        H(q,\dot{q},p_{q}, p_{\dot{q}}) = & \left( \frac{(1 - \cos(2 \theta))^3 (-\frac{1}{8} \dot{\phi} \dot{\theta} p_{\dot{\phi}} \sin(2\theta) + \frac{1}{4} \dot{\phi} \dot{\theta} p_{\dot{\psi}} \sin(\theta) + \frac{1}{4} \dot{\psi} \dot{\theta} p_{\dot{\phi}} \sin(\theta) }{\cos(2 \theta) + 1} \right. \\
        & + \frac{- \frac{1}{8} \dot{\psi} \dot{\theta} p_{\dot{\psi}} \sin(2 \theta) + \frac{1}{4} p_{\dot{\phi}}^2 - \frac{1}{2} p_{\dot{\phi}} p_{\dot{\psi}} \cos(\theta) + \frac{1}{4} p_{\dot{\psi}}^2)}{\cos(2 \theta) + 1} \\
        & - \frac{1}{2} (-p_{\dot{\phi}} + \frac{p_{\dot{\psi}}}{\cos \theta})^2 \sin^4(\theta) - \frac{1}{2} (- \frac{p_{\dot{\phi}}}{\cos\theta} + p_{\dot{\psi}})^2 \sin^4 \theta +  (p_{\dot{\phi}}^2 - p_{\dot{\phi}} p_{\dot{\psi}} \cos(\theta) - \frac{p_{\dot{\phi}} p_{\dot{\psi}}}{\cos\theta} + p_{\dot{\psi}}^2) \sin^4\theta \\
        &  \left.  + (\dot{\phi} p_{\phi} + \dot{\psi} p_{\psi} + \dot{\theta} p_{\theta} - \frac{1}{2} p_{\dot{\theta}}^2 - p_{\dot{\theta}} (\dot{\phi} \dot{\psi} \sin(\theta) - p_{\dot{\theta}})) \sin(\theta)^6 \tan(\theta)^2\right)(\sin(\theta)^6 \tan(\theta)^2)^{-1}
    \end{split}
\end{equation*}
\end{tiny}
where $(q,\dot{q},p_{q}, p_{\dot{q}}) = (\theta, \phi, \psi, \dot{\theta}, \dot{\phi}, \dot{\psi}, p_{\theta}, p_{\phi}, p_{\psi}, p_{\dot{\theta}}, p_{\dot{\phi}}, p_{\dot{\psi}})$ are the coordinates on $T^{*}(T\St{3,2})$.

\subsection{The symplectic numerical methods}

Using formula \eqref{cubic}, one can evaluate points on the geometric cubic polynomial at different values of the parameter $t$ (only requires computing exponentials of skew-symmetric matrices and logarithms of orthogonal matrices), in order to compare the results with curves obtained using different approaches.

 Our objective is to compare the curves obtained in Theorem \ref{theo} with  the approximate Riemannian cubic polynomials obtained using geometric integrators constructed upon Proposition \ref{Prop:symplectic:integrator}.

In the following comparison of numerical methods, we have integrated a higher-order Runge-Kutta method to approximate the Riemannian cubic polynomials in a single chart for the sphere and for the Stiefel manifold $\St{3,2}$ and use it as the benchmark to compare different methods. In the sphere, departing from initial conditions $(\theta_0,\phi_{0})=(\frac{\pi}{2}, \pi)$ with velocity $(\dot{\theta}_0,\dot{\phi}_{0})=(0.1, 0.2)$ and acceleration and jerk equal to $(1, 0.5)$ and $(0.1, 0.2)$, respectively, we integrate Riemannian cubic polynomials for 1 second, split into $N$ evenly spaced time-steps $h>0$. This process generates a discrete flow on $TTQ$, that is, a sequence of positions, velocities, accelerations and jerks denoted by $\{(q_{k}, \dot{q}_{k}, \ddot{q}_{k}, q^{(3)}_{k})\}_{k=0}^{N}$ approximating the Riemannian cubic polynomial. 

\begin{figure}[htb!]
    \centering
    \includegraphics[scale=0.35]{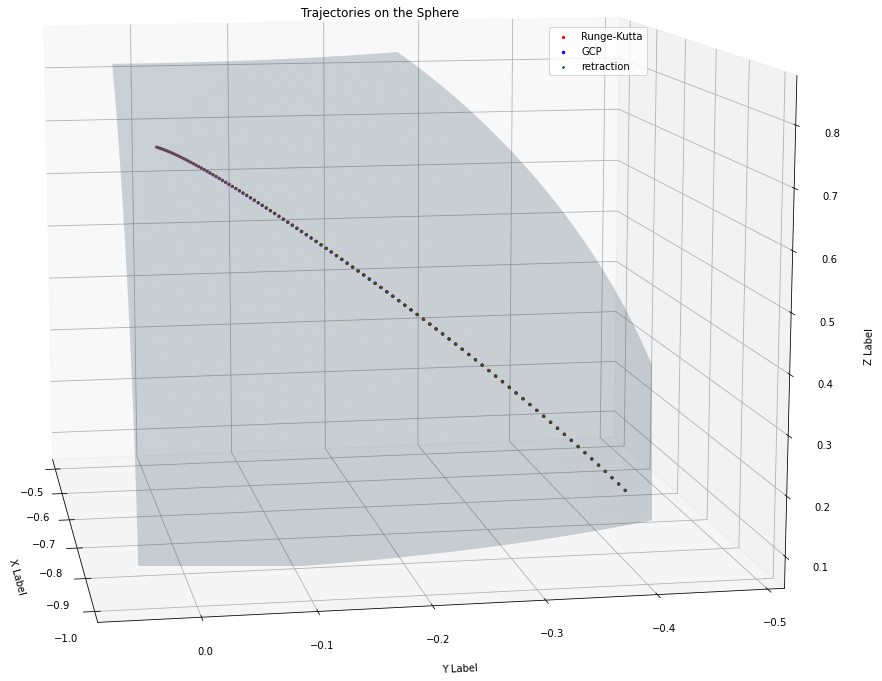}
    \caption{The trajectories on the sphere of the Runge-Kutta method, adjusted de Casteljau algorithm (GCP) and {retraction} based numerical method.}
    \label{trajectories}
\end{figure}

Using this flow, we fix the endpoints $q_{0}$ and $q_{N}$ in the sphere $S^{2}$ as well as the boundary velocities $\dot{q}_{0}$ and $\dot{q}_{N}$. We implemented adjusted de Casteljau algorithm between these points. See Figure \ref{trajectories} {for} the plot of the trajectories.

Then we proceed to compare it with two different retraction based numerical integrators: initial point and mid-point discretization maps from Example \ref{intial:point:discretization} and \ref{mid:point:discretization}, respectively. Both these schemes were implemented in the single chart described in the previous section and in conjunction with a shooting method to enforce the desired boundary conditions. 

The strategy followed for $\St{3,2}$ was very similar to the sphere. Departing from initial conditions $\theta_0 = \pi/4, \phi_0 = \pi, \psi_0 = \pi$, with velocity $\dot{\theta}_{0}=0.1$, $\dot{\phi}_{0}=0.1$ and $\dot{\psi}_{0}=0.05$ and acceleration and jerk equal to $(1, 0.3, 0.5)$ and $(0.1, 0.1, 0.05)$, respectively, we integrated the Riemannian cubic polynomial for 1 second, split into $N$ evenly spaced time-steps $h>0$. 

\subsection{Discussion}

\begin{figure}[h!]
    \centering
    \begin{subfigure}[b]{0.49\textwidth}
        \centering
        \includegraphics[scale=0.43]{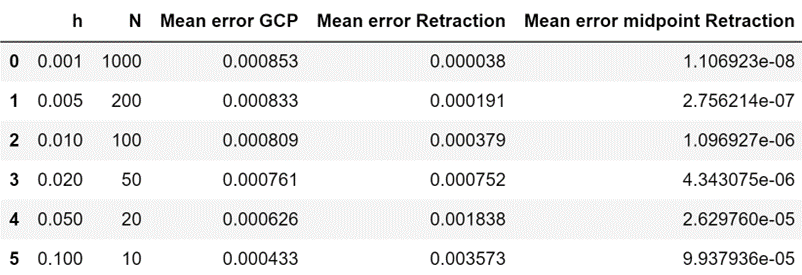}
        \caption{In the sphere.}
        \label{fig:table}
    \end{subfigure}
    \hfill
    \begin{subfigure}[b]{0.49\textwidth}
        \centering
        \includegraphics[scale=0.443]{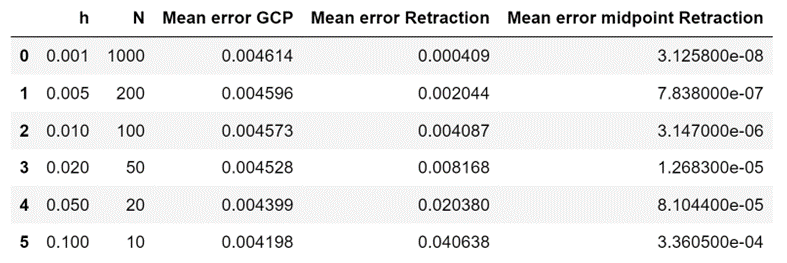}
        \caption{In the Stiefel manifold $St_{3,2}$.}
        \label{fig:table:32}
    \end{subfigure}
    \caption{Table showing the comparison of mean error of adjusted Casteljau algorithm (GCP) with the mean error achieved by a retraction based method using different discretization maps (Retraction and Mid-point Retraction)}
\end{figure}

The first conclusion that we may draw is that the accuracy of the adjusted de Casteljau algorithm on these manifolds is well below that of Casteljau algorithm in Euclidean spaces. Despite this, the adjusted de Casteljau algorithm approximates geometric cubic polynomials reasonably well having relative mean error of around $0.080 \%$ in the sphere and $0.45 \%$ in the Stiefel manifold $St_{3,2}$. The relative error is computed taking the ratio between the average mean error in Figures \ref{fig:table} and \ref{fig:table:32}, respectively, and the maximum distance between two points in these manifolds.

\begin{figure}[h!]
    \centering
    \begin{subfigure}[b]{0.49\textwidth}
        \centering
        \includegraphics[scale=0.49]{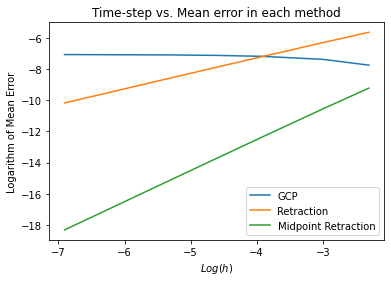}
        \caption{In the sphere.}
        \label{fig:logerror}
    \end{subfigure}
    \hfill
    \begin{subfigure}[b]{0.49\textwidth}
        \centering
        \includegraphics[scale=0.49]{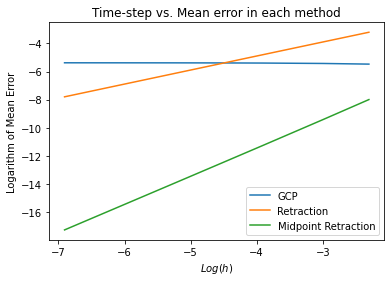}
        \caption{In the Stiefel manifold $St_{3,2}$.}
        \label{fig:logerror:32}
    \end{subfigure}
    \caption{Comparison of mean error of adjusted Casteljau algorithm (GCP) with the mean error achieved by a retraction based method using different discretization maps (Retraction and Mid-point Retraction)}
\end{figure}

On the other hand, the error of retraction-based numerical integrators is an increasing function of the time-step. Ultimately, the mean error of these methods approaches zero at the expense of an increase in computational effort. Still, this property makes retraction-based integrators ideal to simulate dynamics near the initial point of the trajectories, at any desired error order. Notice also, that this is not the case with the adjusted de Casteljau algorithm which is independent of the choice of time-step $h$.

Using retraction maps that are better approximations of the Riemannian exponential map on the configuration manifold, we conjecture that the associated retraction map based integrators will possess a smaller mean error than the mid-point retraction-based symplectic integrator. This trend can already be observed in Figures \ref{fig:logerror} and \ref{fig:logerror:32}: the mid-point discretization map produced a significantly lower error than the initial point discretization.

Finally, retraction-based symplectic integrators require much more computational effort than the adjusted de Casteljau algorithm which runs reasonably fast. Moreover, it is especially suitable for solving boundary value problems while retraction-based methods are initial value problems. Thus, using them to solve boundary value problems requires the use of a shooting method that slows down even further their time cost. 

The retraction map approach faces two challenges: the first one is that, for the time being, they are only suitable to integrate Riemannian cubic polynomials on a single chart. The second is that this approach involves solving several implicit algebraic equations, which makes it very slow and subject to numerical errors. To improve performance, one might use the intrinsic geometric structure of the manifold, in the case that it is either a Lie group or a homogeneous manifold; or by considering it as a submanifold of a higher dimensional Euclidean space, which {will} be analyzed in a future work. In a cost-benefit analysis, the adjusted de Casteljau algorithm revealed to be much faster approximating the Riemannian cubic polynomial satisfying the boundary value problem, though in arbitrary Riemannian manifolds this approximation is far from being as accurate as in the Euclidean case, even in the case of the sphere, where the method relies upon interpolating geodesics, suggesting that one possible direction of research is an alternative choice of control points.

{
\section*{Acknowledgements}
 A. Anahory Simoes and L. Colombo acknowledge the support of the grant PID2022-137909NB-C21 funded by MCIN/AEI/ 10.13039/501100011033. F. Silva Leite is grateful to Funda\c{c}\~ao para a Ci\^encia e Tecnologia (FCT) for the financial support to the project UIDB/00048/2020 (https://doi.org/10.54499/UIDB \\
 /00048/2020). 
  }

\bibliographystyle{plain} 
\bibliography{controlo2024}
%
%

\end{document}